\documentclass[12pt]{article}
\usepackage{amsthm,amsfonts,amssymb,amscd}

\begin{document}


\begin{center}
\large \bf Birational geometry of algebraic varieties,\\
fibred into Fano double spaces
\end{center}\vspace{0.3cm}

\centerline{A.V.Pukhlikov}\vspace{0.3cm}

\parshape=1
3cm 10cm \noindent {\small \quad\quad\quad \quad\quad\quad\quad
\quad\quad\quad {\bf }\newline We develop the quadratic technique
of proving birational rigidity of Fano-Mori fibre spaces over a
higher-dimensional base. As an application, we prove birational
rigidity of generic fibrations into Fano double spaces of
dimension $M\geqslant 4$ and index one over a rationally connected
base of dimension at most $\frac12 (M-2)(M-1)$. An estimate for
the codimension of the subset of hypersurfaces of a given degree
in the projective space with a positive-dimensional singular set
is obtained, which is close to the optimal one.

Bibliography: 15 titles.} \vspace{0.3cm}

\begin{flushleft}
14E05, 14E07
\end{flushleft} \vspace{0.3cm}

\vspace{0.3cm}

\section*{Introduction}

{\bf 0.1. Statement of the main result.} In \cite{Pukh15a}
birational rigidity was shown for two large classes of
higher-dimensional Fano-Mori fibre spaces: generic fibrations into
double spaces of index one and dimension $M \geqslant 5$ when the
dimension of the base does not exceed $\frac12 (M-4)(M-1)-1$ and
generic fibrations into hypersurfaces of index one and dimension
$M-1\geqslant 9$ when the dimension of the base does not exceed
$\frac12(M-7)(M-6)-6$ (in both cases under the assumption of
sufficient twistedness over the base). For Fano-Mori fibre spaces
over the projective line the question of birational rigidity is
studied well enough, see \cite[Chapter 5]{Pukh13a}. However, one
should note that almost all results on birational rigidity of
Fano-Mori fibre spaces over the line were obtained by means of the
{\it quadratic} technique (that is, via analysis of the
singularities of the self-intersection of a mobile linear system,
defining the birational map), whereas the main result of
\cite{Pukh15a} was obtained by means of the {\it linear} technique
(that is, via direct analysis of the singularities of the linear
system itself, without using the quadratic operation of taking the
self-intersection). The quadratic technique requires less
restrictions on the variety underconsideration and for that reason
makes it possible to embrace a considerably large class of
rationally connected varieties. In many respects it is more
efficient (at least, at the present stage of of the theory of
birational rigidity).\vspace{0.1cm}

The aim of the present paper is to develop the quadratic technique
of studying birational geometry of Fano-Mori fibre spaces over a
higher-dimensional base and apply it to fibrations into Fano
double spaces of index one, considerably improving the result of
\cite{Pukh15a} for that class of varieties: we show birational
rigidity of generic fibrations into Fano double spaces of
dimension $M\geqslant 4$ and index one over a rationally connected
base of dimension at most $\frac12 (M-2)(M-1)$. This result,
considerably increasing the admissible dimension of the base of
the fibre space, is obtained by means of the quadratic technique
of counting multiplicities (see \cite[Chapter 5]{Pukh13a}), which
was not used in \cite{Pukh15a}.\vspace{0.1cm}

Let us make the precise statements.\vspace{0.1cm}

We consider a Fano-Mori fibre space $\pi\colon V\to S$, where the
base $S$ is non-singular, the variety $V$ is factorial and has at
most terminal singularities, the antocanonical class $(-K_V)$ is
relatively ample and
$$
\mathop{\rm Pic} V = {\mathbb Z} K_V\oplus \pi^* \mathop{\rm Pic}
S.
$$
We say that a fibre $F=F_s=\pi^{-1}(s)$, $s\in S$, satisfies the
{\it condition} ($h$), if for any irreducible subvariety $Y\subset
F$ of codimension 2 and any point $o\in Y$ the inequality
$$
\frac{\mathop{\rm mult}\nolimits_o Y}{\mathop{\rm deg} Y}\leqslant
\frac{4}{\mathop{\rm deg} F}
$$
holds, where the degrees are understood in the sense of the
anticanonical class, that is,
$$
\mathop{\rm deg} Y=\left(Y\cdot (-K_V)^{\mathop{\rm dim} Y}\right)
$$
and
$$
\mathop{\rm deg} F=\left(F\cdot (-K_V)^{\mathop{\rm dim}
F}\right),
$$
and the {\it condition} ($hd$), if for any mobile linear system
$\Delta\subset |-n(K_V|_F)|$ and any irreducible subvariety
$Y\subset F$ of codimension 2 the inequality
$$
\mathop{\rm mult}\nolimits_Y \Delta\leqslant n
$$
holds. Further, we say that a fibre $F$ satisfies the {\it
condition} ($v$), if for any prime divisor $Y\subset F$ and any
point $o\in F$ of this fibre the inequality
$$
\frac{\mathop{\rm mult}\nolimits_o Y}{\mathop{\rm deg} Y}\leqslant
\frac{2}{\mathop{\rm deg} F}
$$
holds. Finally, we say that the fibre space $V/S$ satisfies the
$K$-{\it condition}, if for any mobile family $\overline{\cal C}$
of curves on the base $S$, sweeping out $S$, and a general curve
$\overline{C}\in \overline{\cal C}$ the class of algebraic cycle
$$
-N (K_V\cdot \pi^{-1}(\overline{C}))-F
$$
of dimension $\mathop{\rm dim} F$ for any $N\geqslant 1$ is not
effective, that is, it is not rationally equivalent to an
effective cycle of dimension $\mathop{\rm dim} F$, and the
$K^2$-{\it condition}, if for any mobile family $\overline{\cal
C}$ of curves on the base $S$, sweeping out $S$, and a general
curve $\overline{C}\in \overline{\cal C}$ the class of algebraic
cycle
$$
N(K_V^2\cdot \pi^{-1}(\overline{C}))-H_F
$$
of dimension $\mathop{\rm dim} F -1$ is not effective for any
$N\geqslant 1$, where $H_F=(-K_V\cdot F)$ is the class of the
anticanonical section of the fibre.\vspace{0.1cm}

The following claim is the main result of the present
paper.\vspace{0.1cm}

{\bf Theorem 0.1.} {\it Assume that $\mathop{\rm dim} F\geqslant
4$ and every fibre $F$ of the projection $\pi$ is a variety with
at most quadratic singularities of rank at least 4, and moreover
$\mathop{\rm codim} (\mathop{\rm Sing} F\subset F)\geqslant 4$.
Assume further that every fibre $F$ satisfies the conditions (h),
(hd) and (v), whereas the fibre space $V/S$ satisfies the
$K$-condition $K^2$-condition.\vspace{0.1cm}

Then the fibre space $V/S$ is birationally rigid: every birational
map $\chi\colon V\dashrightarrow V'$ onto the total space of
rationally connected fibre space $V'/S'$ is fibre-wise, that is,
there is a rational dominant map $\beta\colon S\dashrightarrow S'$
such that the following diagram commutes:}
$$
\begin{array}{rcccl}
   & V & \stackrel{\chi}{\dashrightarrow} & V' & \\
\pi & \downarrow &   &   \downarrow & \pi' \\
   & S & \stackrel{\beta}{\dashrightarrow} & S'.
\end{array}
$$

(Recall that a morphism of projective algebraic varieties
$\pi'\colon V'\to S'$ is a rationally connected fibre space if the
base $S'$ and the general fibre ${\pi'}^{-1}(s')$, $s'\in S'$, are
rationally connected.)\vspace{0.1cm}

Theorem 0.1 implies immediately the following claim.\vspace{0.1cm}

{\bf Corollary 0.1.} {\it In the assumptions of Theorem 0.1 on the
variety $V$ there are no structures of a rationally connected
fibre space over a base of dimension higher than $\mathop{\rm dim}
S$. In particular, the variety $V$ is non-rational. Any birational
self-map of the variety $V$ is fibre-wise and induces a birational
self-map of the base $S$, so that there is a natural homomorphism
of groups $\rho\colon \mathop{\rm Bir} V\to \mathop{\rm Bir} S$,
the kernel of which $\mathop{\rm Ker} \rho$ is the group
$\mathop{\rm Bir} F_{\eta} = \mathop{\rm Bir} (V/S)$ of birational
self-maps of the generic fibre $F_{\eta}$ (over the non-closed
generic point $\eta$ of the base $S$), whereas the group
$\mathop{\rm Bir} V$ is an extension of the normal subgroup
$\mathop{\rm Bir} F_{\eta}$ by the group  $\Gamma=\rho(\mathop{\rm
Bir} V)\subset \mathop{\rm Bir} S$:}
$$
1\to \mathop{\rm Bir} F_{\eta} \to \mathop{\rm Bir} V\to \Gamma
\to 1.
$$

Recall that in \cite{Pukh15a} the following fact was
shown.\vspace{0.1cm}

{\bf Theorem 0.2.} {\it Assume that a Fano-Mori fibre space
$\pi\colon V\to S$ satisfies the following
conditions:\vspace{0.1cm}

{\rm (i)} every fibre $F_s={\pi}^{-1}(s)$, $s\in S$, is a
factorial Fano variety with at most terminal singularities and the
Picard group $\mathop{\rm Pic} F_s = {\mathbb Z} K_{F_s}$, where
$F_s$ has complete intersection singularities and $\mathop{\rm
codim} (\mathop{\rm Sing} F\subset F)\geqslant 4$,\vspace{0.1cm}

{\rm (ii)} for every effective divisor $D\in |-nK_{F_s}|$ on an
arbitrary fibre $F_s$ the pair $(F_s,\frac{1}{n} D)$ is log
canonical, and for any mobile linear system $\Sigma_s\subset
|-nK_{F_s}|$ the pair $(F_s,\frac{1}{n} \Sigma_s)$ is canonical
(that is, the pair $(F_s,\frac{1}{n} D)$ is canonical for a
general divisor $D\in\Sigma_s$),\vspace{0.1cm}

{\rm (iii)} for any mobile family $\overline{\cal C}$ of curves on
the base $S$, sweeping out $S$, and a general curve
$\overline{C}\in \overline{\cal C}$ the class of algebraic cycle
of dimension $\mathop{\rm dim} F$ for any positive $N\geqslant 1$
$$
-N (K_V\cdot \pi^{-1}(\overline{C}))-F
$$
(where $F$ is the fibre of the projection $\pi$) is not effective,
that is, it is not rationally equivalent to an effective cycle of
dimension $\mathop{\rm dim} F$.\vspace{0.1cm}

Then any birational map $\chi\colon V\dashrightarrow V'$ onto the
total space of a rationally connected fibre space $V'/S'$ is
fibre-wise, that is, there is a rational dominant map $\beta\colon
S\dashrightarrow S'$ such that the following diagram commutes:}
$$
\begin{array}{rcccl}
   & V & \stackrel{\chi}{\dashrightarrow} & V' & \\
\pi & \downarrow &   &   \downarrow & \pi' \\
   & S & \stackrel{\beta}{\dashrightarrow} & S'.
\end{array}
$$

Let us compare the assumptions of Theorems 0.1 and 0.2. The
canonicity of the pair $(F_s,\frac{1}{n}D)$ in the condition (ii)
of Theorem 0.2 (that is, essentially the birational superrigidity
of the fibre $F_s$) follows from the conditions ($h$) and ($hd$)
in Theorem 0.1 (and is actually equivalent to them as the main
method of proving birational superrigidity of a primitive Fano
variety is the application of the $4n^2$-inequality combined with
the exclusion of maximal subvarieties of codimension two, see
\cite[Chapter 2]{Pukh13a}). The log canonicity of the pair
$(F_s,\frac{1}{n}D)$ in the condition (ii) of Theorem 0.2 is
replaced in Theorem 0.1 by the condition ($v$), which for certain
classes of Fano varieties is much easier to check. Finally, in
Theorem 0.1 a new global condition for the Fano-Mori fibre space
$V/S$ is added, the $K^2$-condition, which is easy to
check.\vspace{0.1cm}

Theorem 0.1 will be applied to fibrations into double spaces of
index one, when the conditions ($h$) and ($v$) hold automatically
by the equality $\mathop{\rm deg} F=2$.\vspace{0.3cm}


{\bf 0.2. Fibrations into double spaces of index one.} We use the
notations of subsection 0.2 of \cite{Pukh15a}: the symbol
${\mathbb P}$ stands for the projective space ${\mathbb P}^M$,
$M\geqslant 4$, and ${\cal W}={\mathbb P}(H^0({\mathbb P},{\cal
O}_{\mathbb P}(2M)))$ is the space of hypersurfaces of degree $2M$
in ${\mathbb P}$. The following general fact is
true.\vspace{0.1cm}

{\bf Theorem 0.3.} {\it The closed algebraic subset of homogeneous
polynomials $f$ of degree $d$ in $(N+1)$ variables, such that the
hypersurface $\{f=0\}\subset {\mathbb P}^N$ has a singular set of
positive dimension, is of codimension at least $(d-2)N$ in the
space} $H^0({\mathbb P}^N,{\cal O}_{{\mathbb
P}^N}(d))$.\vspace{0.1cm}

{\bf Proof} is given in \S 3.\vspace{0.1cm}

The following theorem is immediately implied by Theorem
0.3.\vspace{0.1cm}

{\bf Theorem 0.4.} {\it There is a Zariski open subset ${\cal
W}_{\rm reg}\subset {\cal W}$, such that any hypersurface $W\in
{\cal W}_{\rm reg}$ has finitely many singular points, each of
which is a quadratic singularity of rank at least 3, and,
moreover, the following estimate holds:}
$$
\mathop{\rm codim} (({\cal W}\setminus {\cal W}_{\rm reg})\subset
{\cal W})\geqslant \frac{(M-2)(M-1)}{2}+1.
$$

{\bf Proof.} Setting in Theorem 0.3 $d=2M$ and $N=M$, we obtain,
that in the complement to a closed subset of codimension $2M(M-1)$
in ${\cal W}$ any hypersurface $W$ has finitely many singular
points. It is easy to check that the closed set of hypersurfaces
$W$ with a quadratic singular point of rank at most 2 or with a
singularity $o\in W$ of multiplicity $\mathop{\rm mult}\nolimits_o
W\geqslant 3$, is of codimension $\frac12(M-2)(M-1)+1$ in the
space ${\cal W}$. This proves the theorem.\vspace{0.1cm}

If $F\to {\mathbb P}$ is a double cover, branched over a
hypersurface $W\in {\cal W}_{\rm reg}$, then $F$ is a factorial
Fano variety with terminal singularities (see \cite{CL},
Subsection 2.1 in \cite{Pukh15a} and Proposition 1.4 below),
satisfying the conditions ($h$) and ($v$) by the equality
$\mathop{\rm deg} F=2$. The condition ($hd$) is easy to show by
the standard methods (see \cite[Chapter 2]{Pukh13a}; for
$M\geqslant 5$ it holds in a trivial way, because for any
irreducible subvariety $Y\subset F$ of codimension 2 the
inequality $\mathop{\rm deg} Y\geqslant 2$ holds). Thus in order
to apply Theorem 0.1, it is sufficient to require every fibre
$F_s$, $s\in S$, to be branched over a regular hypersurface
$W_s\in {\cal W}_{\rm reg}$, and the fibre space $V/S$ to satisfy
the $K$-condition and the $K^2$-condition.\vspace{0.1cm}

In the notations of Subsection 0.2 of \cite{Pukh15a} let $S$ be a
non-singular rationally connected variety of dimension
$\mathop{\rm dim} S\leqslant \frac12 (M-2)(M-1)$. Let ${\cal L}$
be a locally free sheaf of rank $M+1$ on $S$ and $X={\mathbb
P}({\cal L})={\bf Proj}\,
\mathop{\oplus}\limits_{i=0}^{\infty}{\cal L}^{\otimes i}$ the
corresponding ${\mathbb P}^M$-bundle. We may assume that ${\cal
L}$ is generated by its global sections, so that the sheaf ${\cal
O}_{{\mathbb P}({\cal L})}(1)$ is also generated by the global
sections. Let $L\in \mathop{\rm Pic} X$ be the class of that
sheaf, so that
$$
\mathop{\rm Pic} X = {\mathbb Z} L\oplus \pi_X^* \mathop{\rm Pic}
S,
$$
where $\pi_X\colon X\to S$ is the natural projection. Take a
general divisor $U\in |2(ML+\pi^*_X R)|$, where $R\in \mathop{\rm
Pic} S$ is some class. If that system is sufficiently mobile, then
by the assumption about the dimension of the base $S$ and by
Theorem 0.4 we may assume that for any point $s\in S$ the
hypersurface $U_s= U\cap \pi^{-1}_X(s)\in {\cal W}_{\rm reg}$, and
for that reason the double space, branched over $U_s$, satisfies
the conditions of Theorem 0.1. Let $\sigma\colon V\to X$ be the
double cover branched over $U$. Set $\pi=\pi_X\circ\sigma\colon
V\to S$, so that $V$ is a fibration into Fano double spaces of
index one over $S$. Recall that the divisor $U\in |2(ML+\pi^*_X
R)|$ is assumed to be sufficiently general.\vspace{0.1cm}

{\bf Theorem 0.5.} {\it Assume that the variety $V$ is general in
the sense of the construction described above and the divisorial
class $(K_S+R)$ is pseudo-effective. Then for the fibre space
$\pi\colon V\to S$ the claims of Theorem 0.1 and Corollary 0.1 are
true. In particular,
$$
\mathop{\rm Bir} V = \mathop{\rm Aut} V = {\mathbb Z}/2{\mathbb Z}
$$
is the cyclic group of order 2.}\vspace{0.1cm}

{\bf Proof.} Since the class $L$ is numerically effective and
$$
\left( (\sigma^* L)^M\cdot F\right)= \left( (\sigma^*
L)^{M-1}\cdot H_F\right)=2,
$$
it is sufficient to check the inequalities
$$
\left( (\sigma^* L)^M\cdot K_V\cdot
\pi^{-1}(\overline{C})\right)\geqslant 0\quad\mbox{and}\quad
\left( (\sigma^* L)^{M-1}\cdot K^2_V\cdot
\pi^{-1}(\overline{C})\right)\leqslant 0.
$$
As it was noted in Subsection 0.2 in \cite{Pukh15a}, the first of
these inequalities up to a positive factor is the inequality
$((K_S+R)\cdot \overline{C})\geqslant 0$, which holds because the
class $(K_S+R)$ is pseudo-effective and the family of curves
$\overline{\cal C}$ is mobile and sweeps out the base $S$. As for
the second inequality, then elementary computations show that up
to a positive factor it can be written as the inequality
$$
2((K_S+R)\cdot \overline{C})+((\mathop{\rm det} {\cal L})\cdot
\overline{C})\geqslant 0,
$$
which is the more so true because the locally free sheaf ${\cal
L}$ is generated by global sections. Q.E.D. for the theorem.
\vspace{0.1cm}

{\bf Remark 0.1.} For fibrations into double spaces of index one
the $K^2$-condition follows from the $K$-condition. Theorem 0.5
makes Theorem 0.3 of \cite{Pukh15a} stronger in respect of the
genericity conditions which should be satisfied for every fibre of
the fibre space $V/S$: in Theorem 0.5 these conditions are weaker,
and for that reason the set ${\cal W}_{\rm reg}$ is larger. This
makes it possible to prove birational rigidity for fibre spaces
over a base of higher dimension, and in particular, for fibrations
onto four-dimensional double spaces.\vspace{0.3cm}


{\bf 0.3. The structure of the paper.} The present paper is
organized in the following way. \S 1 contains mostly the first
part of the proof of Theorem 0.1: we construct a modification of
the base $S^+\to S$ such that on the pull back $\pi_+\colon V^+\to
S^+$ of the original fibre space onto $S^+$, the centre of each
maximal singularity covers a divisor on $S^+$ (this procedure is
often referred to as {\it flattening} the maximal singularities).
These arguments are similar to the arguments of \S 1 in
\cite{Pukh15a}, however, in contrast to \cite{Pukh15a}, here they
give no proof of the main theorem, but only show the existence of
a {\it supermaximal singularity} (under the assumption that the
claim of Theorem 0.1 does not hold). The latter concept plays an
important role in the proof of birational superrigidity of fibre
spaces over ${\mathbb P}^1$, see \cite[Chapter 5]{Pukh13a}; here
we extend it to the case of fibrations over a base of arbitrary
dimension. We complete \S 1, studying quadratic singularities, the
rank of which is bounded from below (we need this to claim
factoriality and terminality of the modified fibre space
$\pi_+\colon V^+\to S^+$).\vspace{0.1cm}

In \S 2 we complete the proof of Theorem 0.1: we exclude the
supermaximal singularity, the existence of which has been shown in
\S 1, whence the claim of Theorem 0.1 follows immediately. The
excluding is achieved by means of the standard technique of
counting multiplicities (see \cite[Chapter 5]{Pukh13a}), adjusted
to the situation under consideration.\vspace{0.1cm}

In \S 3 we obtain an estimate for the codimension of the closed
set of hypersurfaces of degree $d$ in ${\mathbb P}^N$ with a
singular set of positive dimension, in the space of all
hypersurfaces of degree $d$ in ${\mathbb P}^N$. The estimate is
close to the optimal one. This is a general and quite useful
result, proved by elementary (but non-trivial) methods of
algebraic geometry; as far as the author knows, this estimate was
not known earlier.\vspace{0.3cm}


{\bf 0.4. Historical remarks and acknowledgements.} The history of
the problems connected with birational rigidity of Fano-Mori fibre
spaces over a base of positive dimension, has been reviewed in the
introduction to \cite{Pukh15a} in a detailed enough way, and we
will not consider it here. We note, however, that the Sarkisov
theorem on conic bundles \cite{S80,S82} was proved by the
quadratic method (the self-intersection of the mobile linear
system, defining the birational map, was considered), although for
that class of varieties the quadratic technique of counting
multiplicities is not needed.\vspace{0.1cm}

The problem of birational rigidity for del Pezzo fibrations over a
base of dimension higher than one is entirely open. However, it is
clear that for that class of varieties it is the quadratic
techniques that is needed, although it is possible that a
combination of the linear and the quadratic method will be
successful. In the direction of computing the possible values of
log canonical thresholds on del Pezzo surfaces a lot of work has
been recently done, see \cite{Ch08,Ch14a,Ch14b,ChK}.\vspace{0.1cm}

Finally, let us point out the recent work \cite{Kr}, where by
means of the results of \cite{Pukh05,Pukh09b} (see also
\cite[Chapter 7]{Pukh13a}) the problem of existence of rationally
connected varieties that are non-Fano type varieties, stated in
\cite{CaGo}, was solved.\vspace{0.1cm}

Various technical points, related to the constructions of the
present paper, were discussed by the author in his talks given in
2009-2014 at Steklov Mathematical Institute. The author thanks the
members of Divisions of Algebraic Geometry and of Algebra and
Number Theory for the interest to his work. The author also thanks
his colleagues in the Algebraic Geometry research group at the
University of Liverpool for the creative atmosphere and general
support.


\section{Maximal and supermaximal \\ singularities}

The contents of this section is the first part of the proof of
Theorem 0.1. In Subsection 1.1 we modify the fibre space $V/S$:
this procedure is similar to \S 1 in \cite{Pukh15a}. As a result,
we obtain a new Fano-Mori fibre space $V^+/S^+$, satisfying all
assumptions of Theorem 0.1 and an additional condition: the centre
on $V^+$ of any maximal singularity covers a divisor on $S^+$. In
Subsection 1.2 we consider the self-intersection of the mobile
linear system $\Sigma$, related to the birational map $\chi$, and
show the existence of a supermaximal singularity. In Subsection
1.3 we make the information about quadratic singularities of a
bounded rank more precise.\vspace{0.3cm}

{\bf 1.1. Modification of the fibre space $V/S$.} In the notations
of Theorem 0.1 fix a birational map $\chi\colon V\dashrightarrow
V'$. Repeating the arguments of Subsection 1.1 in \cite{Pukh15a},
consider an arbitrary very ample linear system
$\overline{\Sigma'}$ on $S'$. Let
$\Sigma'=(\pi')^*\overline{\Sigma'}$ be its pull back onto $V'$,
so that the divisors $D'\in\Sigma'$ are composed from the fibres
of the projection $\pi'$, and for that reason for any curve
$C\subset V'$ that is contracted by the projection $\pi'$, we have
$(D'\cdot C)=0$; the linear system $\Sigma'$ is obviously mobile.
Set
$$
\Sigma=(\chi^{-1})_*\Sigma'\subset |-nK_V+\pi^*Y|
$$
to be its strict transform on $V$, where $n\in{\mathbb Z}_+$.
Obviously, the map $\chi$ is fibre-wise if and only if $n=0$.
Therefore, if $n=0$, then the claim of Theorem 0.1 holds. So let
us assume that $n\geq 1$ and show that this assumption leads to a
contradiction.\vspace{0.1cm}

As was shown in \cite[Lemma 1.1]{Pukh15a}, for any mobile family
of curves $\overline{C}\in\overline{\cal C}$ on $S$, sweeping out
$S$, the inequality $(\overline{C}\cdot Y)\geqslant 0$
holds.\vspace{0.1cm}

Following \cite{Pukh15a}, we call a prime divisor $E$ over $V$ a
{\it maximal singularity} of the birational map $\chi$, if its
image on $V'$ is a prime divisor, covering the base $S'$, and the
Noether-Fano inequality holds:
$$
\varepsilon(E)=\mathop{\rm ord}\nolimits_E \Sigma -na(E)>0,
$$
where $a(E)$ is the discrepancy of $E$ with respect to $V$. In
\cite[Proposition 1.1]{Pukh15a} it was shown that maximal
singularities exist. Let ${\cal M}$ be the (finite) set of all
maximal singularities.\vspace{0.1cm}

In the proof of the existence of maximal singularities an
important role is played by a {\it very mobile} family ${\cal C}'$
of rational curves on the variety $V'$. Recall \cite[Subsection
1.1]{Pukh15a}, that a family of rational curves ${\cal C}'$ on
$V'$ is very mobile if the curves $C'\in {\cal C}'$ are contracted
by the projection $\pi'$, sweep out a dense open subset in $V'$,
do not intersect the set of indeterminancy of the map
$\chi^{-1}\colon V'\dashrightarrow V$, and a general curve $C'\in
{\cal C}'$ intersects the image of each maximal singularity $E\in
{\cal M}$ transversally at points of general position. Let us fix
a very mobile family of curves on $V'$. Its strict transform on
$V$ we denote by the symbol ${\cal C}$, and its projection
$\pi({\cal C})$ on $S$ by the symbol $\overline{{\cal C}}$.
Further, the following fact is true.\vspace{0.1cm}

{\bf Proposition 1.1.} {\it For every maximal singularity
$E\subset{\cal M}$ its centre
$$
\mathop{\rm centre}(E,V)=\varphi(E)
$$
on $V$ does not cover the base: $\pi(\mathop{\rm
centre}(E,V))\subset S$ is a proper closed subset of the variety
$S$.}\vspace{0.1cm}

{\bf Proof.} Although the statement of this proposition repeats
the statement of Proposition 1.2 in \cite{Pukh15a} word for word,
a new proof is needed, since the assumptions are different. Again
it is sufficient to show that the restriction $\Sigma|_F$ of the
linear system $\Sigma$ onto a fibre $F=\pi^{-1}(s)$ of general
position has no maximal singularities (in the standard, weaker
sense, see \cite[Chapter 2]{Pukh13a}). This follows immediately
from the conditions ($h$) and ($hd$), which are satisfied for the
variety $V$. Q.E.D. for the proposition.\vspace{0.1cm}

Now let us construct, following \cite[Subsection 1.2]{Pukh15a}, a
modification of the base $\sigma_S\colon S^+\to S$ and the
corresponding modification of the total space
$$
\sigma_S\colon V^+=V\times_SS^+\to V
$$
of the fibre space $V/S$, such that the new fibre space
$\pi_+\colon V^+\to S$ satisfies the following conditions:
\begin{itemize}

\item the base $S^+$ is non-singular,

\item for every singularity $E$ of the birational map
$\chi\circ\sigma\colon V^+\dashrightarrow V'$, which is realized
on $V'$ by a divisor, covering the base $S'$, its centre on $V^+$
covers a divisor on $S^+$, that is,
$$
\mathop{\rm codim}(\pi_+(\mathop{\rm centre}(E,V^+))\subset
S^+)=1.
$$
\end{itemize}
The modification $\sigma_S$ is constructed as a sequence of blow
ups with non-singular centres. Since the fibre of the fibre space
$V^+/S^+$ over a point $p\in S^+$ is naturally isomorphic to the
fibre of the original fibre space $V/S$ over the point
$\sigma_S(p)\in S$, and the base $S^+$ is non-singular, by the
assumption about the singularities of the fibres of the original
fibre space $V/S$, the variety $V^+$ has at most quadratic (in
particular, hypersurface) singularities of rank at least 4, and
moreover, $\mathop{\rm codim}(\mathop{\rm Sing}V^+\subset
V^+)\geqslant 4$, so that the variety $V^+$ is factorial and
terminal. Obviously,
$$
\mathop{\rm Pic}V^+={\mathbb Z}K_+\oplus\pi^*_+\mathop{\rm
Pic}S^+,
$$
so that $V^+/S^+$ is again a Fano-Mori fibre space. Let
$\overline{{\cal T}}$ be the set of all $\sigma_S$-exceptional
prime divisors on $S^+$ and ${\cal T}$ the set of all
$\sigma$-exceptional prime divisors on $V^+$. The map
$$
{\cal T}\ni T\mapsto\pi_+(T)=\overline{T}\ni\overline{\cal T}
$$
is a bijection between ${\cal T}$ and $\overline{\cal T}$, the
inverse map is
$$
\overline{\cal
T}\ni\overline{T}\mapsto\pi^{-1}_+(\overline{T})=T\in{\cal T}.
$$
Obviously, $\mathop{\rm Pic}S^+=\sigma^*_S\mathop{\rm
Pic}S\bigoplus \bigoplus\limits_{\overline{T}\in \overline{\cal
T}}{\mathbb Z}\overline{T}$ and a similar equality is true for
$\mathop{\rm Pic}V^+$.\vspace{0.1cm}

{\bf Proposition 1.2.} {\it For the Fano-Mori fibre space
$V^+/S^+$ the $K$-condition and the $K^2$-condition
hold.}\vspace{0.1cm}

{\bf Proof.} Let $\overline{\cal R}$ be a mobile family of curves
on $S^+$, sweeping out $S^+$, and ${\overline R}\in{\overline{\cal
R}}$ a general curve. Then, obviously, $\sigma_S(\overline{\cal
R})$ is a mobile family of curves on $S$, sweeping out $S$, and
$\sigma_S(\overline{R})$ is a general curve in that family. We
have
$$
K_{S^+}=\sigma^*_SK_S+\sum_{\overline{T}\in \overline{\cal
T}}a_T\overline{T}
$$
and, respectively,
$$
K_+=\sigma^*K_V+\sum_{T\in{\cal T}}a_TT
$$
(the discrepancies of the prime divisors $\overline{T}$ and
$T=\pi^{-1}_+(\overline{T})$ with respect to $S$ and $V$, are
obviously equal), and moreover, $a_T>0$ for all $T\in{\cal T}$.
Let us consider the class of an algebraic cycle
$$
\sigma_*[-N(K_+\cdot\pi^{-1}_+(\overline{R}))-F]=
-N(K_V\cdot\pi^{-1}(\sigma_S(\overline{R})))-\alpha F,
$$
where $\alpha=N\sum\limits_{\overline{T} \in\overline{\cal
T}}a_T(\overline{T}\cdot\overline{R}) +1\geqslant 1$. Since for
the fibre space $V/S$ the $K$-condition is satisfied, we can see
from here that it is satisfied for $V^+/S^+$, too. Let us consider
the $K^2$-condition. Writing out explicitly $K^2_+$, we get:
$(K^2_+\cdot\pi^{-1}_+(\overline{R}))=$
$$
=(\sigma^*K^2_V\cdot\pi^{-1}_+(\overline{R}))+2(\sigma^*K_V\cdot
\left(\sum\limits_{T\in{\cal
T}}a_TT\right)\cdot\pi^{-1}_+(\overline{R}))=
$$
$$=(\sigma^*K^2_V\cdot\pi^{-1}_+(\overline{R}))-
2\sum\limits_{\overline{T}\in\overline{\cal
T}} a_T(\overline{T}\cdot\overline{R})H_F
$$
(since $((\sum\limits_{T\in{\cal T}}a_T
T)^2\cdot\pi^{-1}_+(\overline{R}))=0$). Therefore,
$$
\sigma_*[N(K^2_+\cdot\pi^{-1}(\overline{R}))-H_F]=
N(K^2_V\cdot\pi^{-1}(\sigma_S(\overline{R})))-\beta F,
$$
where $\beta=2N\sum\limits_{\overline{T}\in\overline{\cal
T}}a_T(\overline{T}\cdot\overline{R}) +1\geqslant 1$. Since the
fibre space $V/S$ satisfies the $K^2$-condition, this implies that
$V^+/S^+$ satisfies it as well. Q.E.D. for the
proposition.\vspace{0.1cm}

Since obviously the map $\chi\circ\sigma\colon V^+ \dashrightarrow
V$ is fibre-wise with respect to the projections $\pi_+$, $\pi'$
if and only if the map $\chi$ is fibre-wise, we will prove Theorem
0.1 for the Fano-Mori fibre space $V^+/S^+$. The fibres of that
fibre space by construction are the fibres of the original fibre
space $V/S$, so that for $V^+/S$ all assumptions of Theorem 0.1
are satisfied.\vspace{0.1cm}

From now on, in order to simplify the notations, we assume that
$V^+/S^+$ is the original Fano-Mori fibre space $V/S$, which now
has a new property: every singularity $E$ of the map $\chi$ (which
is still not fibre-wise), the centre of which on $V'$ is
divisorial and covers the base $S'$, has on the variety $V$ the
centre $\mathop{\rm centre}(E,V)$, covering a prime divisor on
$S$. In particular, this is true for every maximal singularity
$E\in{\cal M}$. In order not to make the text more difficult to
read using by new symbols, we will use the symbols ${\cal T}$,
$\overline{\cal T}$ in the new sense: $\overline{\cal T}$ is the
set of such prime divisors $\overline{T}$ on the base $S$, that
for some maximal singularity $E\in{\cal M}$ we have
$\pi(\mathop{\rm centre}(E,V))=\overline{T}$, and ${\cal T}$ is
the set of preimages $T=\pi^{-1}(\overline{T})$ of those divisors
on $V$. The projection $\pi$ gives a one-to-one correspondence
between the sets ${\cal T}$ and $\overline{\cal T}$. Let
$$
\tau\colon{\cal M}\to{\cal T}
$$
be the map, relating to a maximal singularity $E\in{\cal M}$ the
divisor $T\in{\cal T}$, containing its centre $\mathop{\rm
centre}(E,V)$, and ${\overline{\tau}}=\pi\circ\tau\colon{\cal
M}\to\overline{\cal T}$, that is to say,
$\overline{\tau}(E)=\pi(\mathop{\rm centre}(E,V))\subset S$. For
$T\in{\cal T}$ set ${\cal M}_T=\tau^{-1}(T)$, so that
$$
{\cal M}=\bigsqcup_{T\in{\cal T}}{\cal M}_T.
$$

{\bf Remark 1.1.} In the situation considered in \cite{Pukh15a},
the modification of the base completes the proof of birational
rigidity of the fibre space, since by the assumption about the
global log canonical threshold of every fibre, no maximal
singularity, the centre of which covers a divisor on the base, can
exist. In this paper the assumption about the global log canonical
threshold is missing, and for that reason the main part of the
proof of Theorem 0.1 starts when the base is modified and the
centre of every maximal singularity covers a divisor on the base.
In the next subsection we carry out some preparatory work for the
subsequent exclusion of maximal singularities.\vspace{0.3cm}


{\bf 1.2. Supermaximal singularities.} For any maximal singularity
$E\in{\cal M}$ set
$$
t_E=\mathop{\rm ord}\nolimits_E T,
$$
where $T=\tau(E)$. By construction, $T\supset\mathop{\rm
centre}(E,V)$, so that $t_E\geqslant 1$. Let
$\varphi\colon\widetilde{V}\to V$ be a birational morphism,
resolving the singularities of the map $\chi$. Every maximal
singularity $E\in{\cal M}$ is realized on the variety
$\widetilde{V}$ by a prime divisor, which we will denote by the
same symbol $E$. By the definition of the numbers $t_E$ we get:
the divisor
$$
\varphi^*T-\sum_{E\in{\cal M}_T}t_EE
$$
is effective and contains none of the maximal singularities
$E\in{\cal M}$ as a component. Now let us consider the strict
transform $\widetilde{\cal C}$ on $\widetilde{V}$ of the mobile
family of curves ${\cal C}$, which was fixed in Subsection 1.1.
For $\widetilde{C}\in\widetilde{\cal C}$ we have:
$$
\left(\left(\varphi^*T-\sum_{E\in{\cal
M}_T}t_EE\right)\cdot\widetilde{C}\right)\geqslant 0.
$$
Set $\nu_E=\mathop{\rm ord}_E\Sigma$ and let $a_E\geqslant 1$ be
the discrepancy of $E$ with respect to $V$. By the symbol
$\widetilde{K}$ we denote the canonical class $K_{\widetilde{V}}$,
so that for the strict transform $\widetilde{\Sigma}$ of the
linear system $\Sigma$ on $\widetilde{V}$ we have
$$
\widetilde{\Sigma}\subset|-n\widetilde{K}+\widetilde{Y}+\Xi|,
$$
where $\widetilde{Y}=\varphi^*\pi^*Y-\sum\limits_{E\in{\cal
M}}\varepsilon(E)E$ (recall that $\varepsilon(E)=\nu_E-na_E$) and
$\Xi$ is a linear combination of $\varphi$-exceptional divisors
$E'$, such that either the centre of $E'$ on $V'$ is a subvariety
of codimension at least 2, or $\varepsilon(E')\leqslant 0$, and
for that reason $(\widetilde{C}\cdot\Xi)\geqslant 0$. Therefore,
the following inequality holds:
\begin{equation}\label{31.10.2015.1}
\sum_{E\in{\cal
M}}\varepsilon(E)(E\cdot\widetilde{C})>(\overline{Y}\cdot\overline{C}).
\end{equation}
Recall that by the $K$-condition $(Y\cdot\overline{C})\geqslant0$.
On the other hand, as we could see a bit earlier, the estimate
\begin{equation}\label{31.10.2015.2}
(\overline{T}\cdot\overline{C})=(\varphi^*T\cdot\widetilde{C})
\geqslant\sum_{E\in{\cal M}_T}t_E(E\cdot\widetilde{C})
\end{equation}
holds.\vspace{0.1cm}

Now let us consider the self-intersection $Z=(D_1\circ D_2)$ of
the mobile linear system $\Sigma$ (where $D_1,D_2\in\Sigma$ are
general divisor which do not have common components due to the
mobility). Let us write this effective algebraic cycle of
codimension 2 in the following way:
$$
Z=Z^h+Z^v+Z^{\emptyset},
$$
where in the sub-cycle $Z^h$ are collected all components $Z$,
covering the base ({\it the horizontal part} of $Z$), in the
sub-cycle $Z^v$ are collected all components of the cycle $Z$ that
are contained in the divisors $T\in{\cal T}$ and cover
$\overline{T}$ ({\it the vertical part} of $Z$), and in the
sub-cycle $Z^{\emptyset}$ are collected all the other components
of the cycle $Z$ (and that part of the cycle  $Z$ is inessential
for us). Obviously, we have the presentation
$$
Z^v=\sum_{T\in{\cal T}}Z^v_T,
$$
where $Z^v_T$ consists of those components of the vertical part,
which are contained in the divisor $T$ and cover
$\overline{T}$.\vspace{0.1cm}

Let $F=F_s=\pi^{-1}(s)$ be the fibre over a point of general
position $s\in\overline{T}$. Since $\mathop{\rm Pic}F={\mathbb
Z}K_F$, we have
$$
Z^v_T|_F\sim-\lambda_TK_F
$$
for some $\lambda_T\in{\mathbb Z}_+$. Therefore,
$$
(Z^v_T\cdot\pi^{-1}(\overline{C}))=\lambda_T(\overline{T}\cdot\overline{C})H_F.
$$

{\bf Definition 1.1.} A maximal singularity $E\in{\cal M}_T$ is
said to be {\it supermaximal}, if the inequality
\begin{equation}\label{02.11.2015.1}
2n\varepsilon(E)>\lambda_T\mathop{\rm ord}\nolimits_ET
\end{equation}
holds.\vspace{0.1cm}

This definition is modelled on the definition of a supermaximal
singularity for Fano fibre spaces over ${\mathbb P}^1$, see
\cite[Chapter 5]{Pukh13a}, and plays the same role.\vspace{0.1cm}

{\bf Proposition 1.3.} {\it A supermaximal singularity
exists.}\vspace{0.1cm}

{\bf Proof.} Since
$$
Z\sim n^2K^2_V+2n((-K_V)\cdot\pi^*Y)+\pi^*(Y^2),
$$
we have
$$
(Z\cdot\pi^{-1}(\overline{C}))=
n^2(K^2_V\cdot\pi^{-1}(\overline{C}))+2n(Y\cdot\overline{C})H_F,
$$
as obviously $(\pi^*(Y^2)\cdot\pi^{-1}(\overline{C}))=0$. On the
other hand,
$$
(Z\cdot\pi^{-1}(\overline{C}))=(Z^h\cdot\pi^{-1}(\overline{C}))
+\left(\sum_{T\in{\cal T}}\lambda_T(\overline{T}
\cdot\overline{C})\right)H_F+\lambda_{\emptyset}H_F
$$
for some $\lambda_{\emptyset}\in{\mathbb Z}_+$. By the
$K^2$-condition we get the inequality
\begin{equation}\label{02.11.2015.2}
2n(Y\cdot\overline{C})\geqslant\sum_{T\in{\cal
T}}\lambda_T(\overline{T}\cdot\overline{C}).
\end{equation}
Combining the inequalities (\ref{31.10.2015.1}),
(\ref{31.10.2015.2}) and (\ref{02.11.2015.2}), we get
$$
2n\sum_{E\in{\cal M}}\varepsilon(E)(E\cdot\widetilde{C})>
\sum_{T\in{\cal T}}\lambda_T\left(\sum_{E\in{\cal
M_T}}t_E(E\cdot\widetilde{C})\right).
$$
Taking into account that the set of maximal singularities ${\cal
M}$ is a disjoint union of the subsets $M_T$, $T\in{\cal T}$, we
see that in the last inequality every maximal singularity appears
only once. Therefore, for some singularity $E\in{\cal M}_T$ the
inequality
$$
2n\varepsilon(E)(E\cdot\widetilde{C})>\lambda_Tt_E(E\cdot\widetilde{C})
$$
holds. Since $(E\cdot\widetilde{C})>0$ for all $E\in{\cal M}$,
this implies the inequality (\ref{02.11.2015.1}). Q.E.D. forthe
proposition.\vspace{0.3cm}

{\bf 1.3. A remark on quadratic singularities.} In \cite[Theorem
4]{EP} and \cite[п. 2.1]{Pukh15a} it was shown that the quadratic
singularities of rank at least $r\geqslant 1$ are stable with
respect to blow ups. This fact can be made more precise in the
following way.\vspace{0.1cm}

{\bf Proposition 1.4.} {\it Assume that an algebraic variety $X$
has at most quadratic singularities of rank at least $r$, and
moreover, the inequality
$$
\mathop{\rm codim}(\mathop{\rm Sing}X\subset X)\geqslant r
$$
holds. Then for any irreducible subvariety $B\subset X$ there is a
Zariski open subset $U\subset X$, such that $U\cap B\neq
\emptyset$ and the blow up $\widetilde{U}\to U$ along the
subvariety $B_U=B\cap U$ has at most quadratic singularities of
rank at least $r$, and the following inequality holds:}
\begin{equation}\label{03.11.2015.1}
\mathop{\rm codim}(\mathop{\rm Sing}\widetilde{U}\subset
\widetilde{U})\geqslant r.
\end{equation}

{\bf Remark 1.2.} In \cite{EP,Pukh15a} the following obvious fact
was used:if a variety $X$ has at most quadratic singularities of
rank at least $r$, then the inequality $\mathop{\rm
codim}(\mathop{\rm Sing}X\subset X)\geqslant r-1$ holds.
Therefore, the codimension of the singular set $\mathop{\rm
Sing}\widetilde{U}$ is at least $r-1$. The proposition stated
above makes the results of \cite{EP,Pukh15a} more precise: the
property of the singular set of the variety $X$ to have
codimension at least $r$ is also stable with respect to blow
ups.\vspace{0.1cm}

{\bf Proof of Proposition 1.4.} By \cite[Theorem 4]{EP} and
\cite[Subsection 2.1]{Pukh15a} we only need to show the inequality
(\ref{03.11.2015.1}). Obviously, we may assume that
$B\subset\mathop{\rm Sing}X$. Arguing as in Subsection 2.1 of
\cite{Pukh15a}, consider a Zariski open subset $U\subset X$, such
that $B_U$ is a non-singular subvariety, and moreover the rank of
quadratic points $b\in B_U$ is constant and equal to $r_1\geqslant
r$. Let $E_U\subset\widetilde{U}$ be the exceptional divisor of
the blow up $\varphi_B\colon\widetilde{U}\to U$ of the subvariety
$B_U$. Obviously, $\varphi_B|_{E_U}\colon E_U\to B_U$ is a
fibration into quadrics of rank $r_1$. It is clear that the set of
singular points $\mathop{\rm Sing}(\widetilde{U}\backslash E_U)$
is of codimension at least $r$. However, a quadric of rank $r_1$
has a singular set of codimension $r_1-1$. Therefore,
$$
\mathop{\rm codim}(\mathop{\rm Sing}E_U\subset E_U)=r_1-1
$$
so that the more so
$$
\mathop{\rm codim}((\mathop{\rm Sing}\widetilde{U}\cap E_U)\subset
E_U)\geqslant r_1-1
$$
and for that reason
$$
\mathop{\rm codim}((\mathop{\rm Sing}\widetilde{U}\cap E_U)\subset
U)\geqslant r_1\geqslant r.
$$
Q.E.D. for the proposition.


\section{Exclusion of supermaximal singularities}

In this section we complete the proof of Thorem 0.1: we show that
a maximal singularity can not exist. For that purpose, we use the
technique of counting multiplicities (Subsection 2.1) in a
modified form, adjusted to varieties with quadratic singularities.
We prove that the multiplicities of the self-intersection of the
mobile linear system $\Sigma$ along the centres of the
supermaximal singularity satisfy a certain quadratic inequality,
which is impossible, as our computations in Subsection 2.2 show.
This contradiction completes the proof of Theorem 0.1. In
Subsection 2.3 we correct a small issue in
\cite{EP}.\vspace{0.3cm}

{\bf 2.1. The technique of counting multiplicities.} Let us fix a
supermaximal singularity $E$ and the corresponding divisor
$T=\pi^{-1}(\overline{T}$). To simplify the notations, we write
$Z^v$ instead of $Z^v_T$ and $\lambda$ instead of $\lambda_T$: the
other singularities and divisors $T'\in{\cal T}$ take no part in
the subsequent arguments. Let
$$
V_K\rightarrow\dots\rightarrow V_i\rightarrow V_{i-1}
\rightarrow\dots\rightarrow V_0=V
$$
be the resolution of the singularity $E$, that is, the sequence of
blow ups $\varphi_{i,i-1}\colon V_i\to V_{i-1}$ of irreducible
subvarieties $B_{i-1}=\mathop{\rm centre} (E,V_{i-1})$ with
exceptional divisors $E_i=\varphi^{-1}_{i,i-1}(B_{i-1})$, where
the last exceptional divisor $E_K$ is the supermaximal singularity
$E$. The set of indices $I=\{1,\dots,K\}$, parameterizing the blow
ups, is the disjoint union
$$
I=I_0\sqcup I_1\sqcup\dots\sqcup I_{M-1},
$$
where $M=\mathop{\rm dim}F$ is the dimension of the fibre and
$i\in I_k$ if and only if $\mathop{\rm dim}B_{i-1}=\mathop{\rm
dim}S-1+k$ (obviously, $\pi\circ\varphi_{i-1,0}(B_{i-1})=T$, so
that $\mathop{\rm dim}B_{i-1}\geqslant\mathop{\rm dim}T$; here
$$
\varphi_{i,j}=\varphi_{j+1,j}\circ\dots\circ\varphi_{i,i-1}\colon
V_i\to V_j
$$
is a composition of elementary blow ups). Certain sets $I_k$ can
be empty. By Proposition 1.4 for $j\in I_{M-2}\cup I_{M-1}$ we
have $B_{j-1}\not\subset\mathop{\rm Sing}V_{j-1}$. For $j\in I$
set
$$
\mu_j=\mathop{\rm mult}\nolimits_{B_{j-1}}V_{j-1}\in\{1,2\},
$$
so that for $j\in I_{M-2}\cup I_{M-1}$ we have $\mu_j=1$. The
strict transform of a subvariety, an effective divisor or a linear
system on $V_j$ we denote by adding the upper index $j$. For a
general divisor $D\in\Sigma$ write
$$
D^j=\varphi^*_{j,j-1}(D^{j-1})-\nu_jE_j.
$$
Let $Z=(D_1\circ D_2)$ be the self-intersection of the mobile
system $\Sigma$. Writing in the usual way (see \cite[Chapter
2]{Pukh13a})
$$
(D^i_1\circ D^i_2)=(D^{i-1}_1\circ D^{i-1}_2)^i+Z_i,
$$
where $Z_i$ is an effective cycle of codimension 2 with the
support inside the exceptional divisor $E_i$, we define the {\it
degree} $d_i$ of the cycle $Z_i$ in the following way. If
$B_{i-1}\not\subset\mathop{\rm Sing}V_{i-1}$, then for a point
$p\in B_{i-1}$ of general position $\varphi^{-1}_{i,i-1}(p)$ is
the projective space ${\mathbb P}^{\delta_i}$ and
$$
d_i=\mathop{\rm deg}(Z_i|_{\varphi^{-1}_{i,i-1}(p)})
$$
is the degree of an effective divisor in that projective space. If
$B_{i-1}\subset\mathop{\rm Sing}V_{i-1}$, then for a general point
$p\in B_{i-1}$ the fibre $\varphi^{-1}_{i,i-1}(p)$ is an
irreducible quadric in the projective space ${\mathbb
P}^{\delta_i+2}$ and $d_i$ is the degree of the effective cycle
$Z_i|_{\varphi^{-1}_{i,i-1}(p)}$ in that projective space. In both
cases $\delta_i$ means the elementary discrepancy $\mathop{\rm
codim}(B_{i-1}\subset V_{i-1})-\mu_i$. As usual, we break the set
$I$ into the {\it lower part}
$$
I_l=I_0\sqcup\dots\sqcup I_{M-2}
$$
and the {\it upper part} $I_u=I_{M-1}$, and set
$$
L=\mathop{\rm max}\{i\in I_l\}.
$$
Finally for $0\leqslant i< j\leqslant L$ set
$$
m_{i,j}=\mathop{\rm mult}\nolimits_{B_{j-1}}(Z^{j-1}_i),
$$
for $i=0$ we write simply $m_j$. The technique of counting
multiplicities (\cite[Chapter 2]{Pukh13a}) gives the system of
equalities
$$
\begin{array}{ccl}
\mu_1\nu^2_1+d_1 & = & m_1,\\
\mu_2\nu^2_2+d_2 & = & m_2+m_{1,2},\\
& \dots &\\
\mu_i\nu^2_i+d_i & = & m_i+m_{1,i}+\dots+m_{i-1,i},\\
& \dots &\\
\mu_L\nu^2_L+d_L & = & m_L+m_{1,L}+\dots+m_{L-1,L}.\\.
\end{array}
$$
Besides, we have the estimate
$$
d_L\geqslant\sum^K_{i=L+1}\nu^2_i.
$$
Let $\Gamma$ be the oriented graph of the resolution of the
singularity $E$, that is, the graph with the set of vertices $I$
and an oriented edge (arrow) joins the vertices $i$ and $j$
(notation: $i\to j$) if and only if $i>j$ and $B_{i-1}\subset
E^{i-1}_j$. Recall \cite[Chapter 2, Definition 2.1]{Pukh13a}, that
a function $a\colon I_l\to{\mathbb R}_+$ is {\it compatible with
the structure of the graph} $\Gamma$, if
$$
a(i)\geqslant\sum_{I_l\ni j\to i}a(j)
$$
for every $i\in I_l$.\vspace{0.1cm}

{\bf Proposition 2.1.} {\it The function}
$$
r_i=r(i)=\mathop{\rm ord}\nolimits_E\varphi^*_{K,i}E_i
$$
{\it is compatible with the structure of the graph
$\Gamma$.}\vspace{0.1cm}

{\bf Proof.} The cartier divisor
$$
\varphi^*_{K,i}E_i-\sum_{I_l\ni j\to i}\varphi^*_{K,j}E_j
$$
is effective, which immediately implies the claim of the
proposition. Q.E.D.\vspace{0.1cm}

Now \cite[Chapter 2, Proposition 2.4]{Pukh13a} gives the
inequality
$$
\sum^L_{i=1}r_im_i\geqslant
\sum^L_{i=1}r_i\mu_i\nu^2_i+r_L\sum^K_{i=L+1}\nu^2_i.
$$
Extending the definition of the numbers $r_i$ to $i\in I_{M-2}$
and using the obvious fact that $r_i$ is non-increasing as a
function of $i$, we get finally:
\begin{equation}\label{11.11.2015.1}
\sum^L_{i=1}r_im_i\geqslant\sum^K_{i=1}r_i\mu_i\nu^2_i.
\end{equation}

{\bf Remark 2.1.} Let $p_{ai}$ be the number of paths in the
oriented graph $\Gamma$ from the vertex $a$ to the vertex $i$ for
$a\neq i$ (so that $p_{ai}=0$ for $a<i$); set $p_{ii}=1$ for all
$i\in I$. Usually (see \cite[Chapter 2]{Pukh13a}) the technique of
counting multiplicities makes use of the numbers $p_{Ki}$ istead
of $r_i$ in the inequalities of the type (\ref{11.11.2015.1}), and
it is easy to see that for $\mu_1=1$ the equality $r_i=p_{Ki}$
holds. If $\mu_1=2$, then $r_1\geqslant p_{K1}$ (see below). The
inequality (\ref{11.11.2015.1}) remains true, if we replace $r_i$
by $p_{Ki}$, however such a modification is hard to use, since it
is the coefficients $r_i$ that appear both in the explicit form of
the Noether-Fano inequality, and in the explicit expression for
$\mathop{\rm ord}_E\varphi^*_{K,0}T$.\vspace{0.1cm}

Set $L_{\rm sing}=\mathop{\rm max}\{1\leqslant i\leqslant
L\,|\,\mu_i=2\}$.\vspace{0.1cm}

{\bf Proposition 2.2.} (i) {\it For $i\geqslant 1+L_{\rm sing}$
the equality $r_i=p_{Ki}$ holds.\, {\rm (ii)} For $1\leqslant
i\leqslant L_{\rm sing}$ the inequality
$$
p_{Ki}\leqslant r_i\leqslant 2p_{Ki}
$$
holds.}\vspace{0.1cm}

{\bf Proof.} The claim (i) is obvious, since for $i\geqslant
1+L_{\rm sing}$ the exceptional divisor $E_i$ is non-singular over
a general point of the subvariety $B_{i-1}$, so that
$$
\mathop{\rm ord}\nolimits_E\varphi^*_{K,i}E_i=\sum_{j\to i}\quad
\mathop{\rm ord}\nolimits_E\varphi^*_{K,j}E_j
$$
and the decreasing induction gives the equality $r_i=p_{Ki}$. For
$i\leqslant L_{\rm sing}$ the fibre of the exceptional divisor
$E_i$ over a point of general position on $B_{i-1}$ is a quadric
of rank at least 4. If for $j\leqslant L_{\rm sing}$, $j>i$, we
have $j\to i$, then, obviously,
\begin{equation}\label{14.11.2015.1}
\varphi^*_{j,j-1}(E^{j-1}_i)=E^j_i+E_j,
\end{equation}
as in the non-singular case. If $j\to i$ for some $j\geqslant
1+L_{\rm sing}$, then two cases are possible:\vspace{0.1cm}

1) $B_{j-1}\not\subset\mathop{\rm Sing}E^{j-1}_i$, and then again
the equality (\ref{14.11.2015.1}) holds,\vspace{0.1cm}

2) $B_{j-1}\subset\mathop{\rm Sing}E^{j-1}_i$, and then the
equality
\begin{equation}\label{14.11.2015.2}
\varphi^*_{j,j-1}(E^{j-1}_i)=E^j_i+2E_j
\end{equation}
holds.\vspace{0.1cm}

We emphasize that if the equality (\ref{14.11.2015.2}) holds, then
$j>L_{\rm sing}$, so that
$$
\mathop{\rm ord}\nolimits_E\varphi^*_{K,j}E_j=p_{Kj}.
$$
For that reason, every path in the graph $\Gamma$ from the top
vertex $K$ to the vertex $i$ gives an input into the number $r_i$,
which is equal to 1 or 2, and the latter takes place if and only
if the path is of the form
$$
i=j_0\leftarrow j_1\leftarrow\dots\leftarrow j_k\leftarrow j_{k+1}
\leftarrow\dots\leftarrow j_m=K,
$$
where $j_k\leqslant L_{\rm sing}$, $j_{k+1}>L_{\rm sing}$ and for
the arrow $j_{k+1}\to j_k$ the case 2), described above, is
realized. Q.E.D. for the proposition.\vspace{0.3cm}


{\bf 2.2. End of the proof of Theorem 0.1.} Recall that above we
defined the elementary discrepancies $\delta_i=\mathop{\rm
codim}(B_{i-1}\subset V_{i-1})-\mu_i$ for $i=1,\dots,K$. Set
$$
L_{\rm fibre}=\mathop{\rm max}\{1\leqslant i\leqslant K \,|\,
B_{i-1}\subset T^{i-1}\}.
$$
For $1\leqslant i\leqslant L_{\rm fibre}$ we define the numbers
$\gamma_i\in{\mathbb Z}$ by the equalities
$$
\varphi^*_{i,i-1}(T^{i-1})=T^i+\gamma_iE_i,
$$
so that $\gamma_i\in\{1,2\}$.\vspace{0.1cm}

{\bf Proposition 2.3.} {\it The following equalities hold:}

(i) {\it the multiplicity of the linear system $\Sigma$ with
respect to $E$ satisfies the relation}
\begin{equation}\label{16.11.2015.1}
\mathop{\rm ord}\nolimits_E\Sigma=\sum^K_{i=1}r_i\nu_i,
\end{equation}

(ii) {\it the multiplicity of the divisor $T$ with respect to $E$
satisfies the relation}
\begin{equation}\label{16.11.2015.2}
\mathop{\rm ord}\nolimits_ET=\sum^K_{i=1}r_i\gamma_i,
\end{equation}

(iii) {\it the discrepancy of $E$ satisfies the relation}
\begin{equation}\label{16.11.2015.3}
a(E)=\sum^K_{i=1}r_i\delta_i.
\end{equation}

{\bf Proof} repeats the arguments in the non-singular case (see
\cite[Chapter 2]{Pukh13a}) word for word, just the number of paths
$p_{Ki}$ should be replaced by the new coefficients $r_i$. We will
show the equality (\ref{16.11.2015.1}); in the other cases the
arguments are similar. We use the induction on $K\geqslant 1$. If
$K=1$, then the equality (\ref{16.11.2015.1}) is obvious. Let
$K\geqslant 2$. For a general divisor $D\in\Sigma$ write:
$$
\varphi^*_{1,0}D=D^1+\nu_1E_1,
$$
so that
$\varphi^*_{K,0}D=\varphi^*_{K,1}D^1+\nu_1\varphi^*_{K,1}E_1$ and
for that reason
$$
\mathop{\rm ord}\nolimits_E\Sigma=\mathop{\rm ord}\nolimits_ED=
\mathop{\rm ord}\nolimits_ED^1+r_1\nu_1.
$$
For $D^1$ the claim of the proposition holds by the induction
hypothesis. The proof is complete.\vspace{0.1cm}

Set $L^*=\mathop{\rm min}(L,L_{\rm fibre})$ and
$$
m^h_i=\mathop{\rm mult}\nolimits_{B_{i-1}}(Z^h)^{i-1}
$$
for $i=1,\dots,L$, and
$$
m^v_i=\mathop{\rm mult}\nolimits_{B_{i-1}}(Z^v)^{i-1}
$$
for $i=1,\dots,L^*$. Now the left-hand side of the inequality
(\ref{11.11.2015.1}) rewrites in the form
\begin{equation}\label{17.11.2015.1}
\sum^L_{i=1}r_im^h_i+\sum^{L^*}_{i=1}r_im^v_i.
\end{equation}
The first component in this sum does not exceed
$$
4n^2\sum^L_{i=1}r_i,
$$
since the sequence of multiplicities $m^h_i$ is not increasing,
and
$$
m^h_i=\mathop{\rm mult}\nolimits_{B_0}Z^h\leqslant \mathop{\rm
mult}\nolimits_{B_0}(Z^h\circ T) \leqslant 4n^2
$$
by the condition $(h)$.  The ``vertical'' component in the sum
(\ref{17.11.2015.1}) by the condition $(v)$ does not exceed the
number
$$
2\lambda\sum^{L^*}_{i=1}r_i\leqslant 2\lambda\mathop{\rm
ord}\nolimits_ET
$$
(see the equality (\ref{16.11.2015.2})), and the right hand side
of the last inequality is strictly smaller than $4ne$, where
$e=\varepsilon(E)$, by the definition of a supermaximal
singularity (the inequality (\ref{02.11.2015.1})). Combining these
estimates, we get that the left hand side of the inequality
(\ref{11.11.2015.1}) is strictly smaller than the expression
$$
4n^2\sum^L_{i=1}r_i+4ne.
$$
Let us consider now the right hand side of the inequality
(\ref{11.11.2015.1}). By the definition of the number
$\varepsilon(E)$ we have:
\begin{equation}\label{18.11.2015.1}
\sum^K_{i=1}r_i\nu_i=n\sum^K_{i=1}r_i\delta_i+e
\end{equation}
(so that in these notations the Noether-Fano inequality takes the
form of the estimate $e>0$). Using the standard methods, it is
easy to check that the minimum of the right hand side of the
inequality (\ref{11.11.2015.1}) on the hyperplane in the space
${\mathbb R}^K_{(\nu_1,\dots,\nu_K)}$, given by the equation
(\ref{18.11.2015.1}), is attained for $\nu_i=\theta/\mu_i$, where
$\theta$ can be found from the equation (\ref{18.11.2015.1}). We
introduce the following notations:
$$
\Sigma_l=\sum^L_{i=1}r_i,\quad \Sigma_u=\sum^K_{i=L+1}r_i,\quad
\Sigma_{\rm sing}= \sum^{L_{\rm sing}}_{i=1}r_i,\quad \Sigma_{\rm
non-sing}=\sum^K_{i=L_{\rm sing}+1}r_i.
$$
In these notations the inequality (\ref{11.11.2015.1}) implies the
estimate
$$
4n^2\Sigma_l+4ne>2\,\frac{(2n\Sigma_l+n\Sigma_u+e)^2} {\Sigma_{\rm
sing}+2\Sigma_{\rm non-sing}}.
$$
Taking into account that $\Sigma_{\rm sing}+\Sigma_{\rm non-sing}=
\Sigma_l+\Sigma_u$, after easy computations we get:
$$
2n^2\Sigma_l\Sigma_{\rm non-sing}+2ne\Sigma_{\rm non-sing}>
2n^2\Sigma^2_l+2n^2\Sigma_l\Sigma_u+n^2\Sigma^2_u+
2ne\Sigma_l+e^2.
$$
However, $\Sigma_{\rm non-sing}\leqslant\Sigma_l+\Sigma_u$, so
that the previous inequality implies the estimate
$$
2ne\Sigma_u>n^2\Sigma^2_u+e^2,
$$
which can not be true. This contradiction excludes the
supermaximal singularity and completes the proof of Theorem
0.1.\vspace{0.3cm}


{\bf 2.3. Birationally rigid Fano hypersurfaces.} In the context
of the proceedings performed in this subsection, let us consider
the problem of estimating the codimension of the set of non-rigid
hypersurfaces of degree $M$ in ${\mathbb P}^M$, which was set and
solved in \cite{EP}. Working on the present paper, the author
detected an incorrectness in that paper in the proof of the
$4n^2$-inequality for Fano hypersurfaces with quadratic
singularities of rank at least 5 for $M\geqslant 5$ (\cite[Section
3]{EP}). In this subsection we explain what was incorrect and how
it should be corrected. Note that the main claim of
(\cite[Proposition 1]{EP}), and the method of its proof are
valid.\vspace{0.1cm}

Recall that in \cite[Section 3]{EP} the following local fact was
shown. Let $X$ be an algebraic variety with quadratic (in
particular, hypersurface) singularities of rank at least 5 (so
that the set of singular points $\mathop{\rm Sing}X$ is of
codimension at least 4 and for that reason the variety $X$ is
factorial), $B\subset\mathop{\rm Sing}X$ an irreducible
subvariety, $\Sigma$ a mobile linear system on $X$, and moreover,
for some $n\geqslant 1$ the pair $(X,\frac{1}{n}\Sigma)$ is not
canonical; more precisely, it has a non canonical singularity $E$
with the centre at $B$. Then the self-intersection $Z=(D_1\circ
D_2)$, $D_i\in\Sigma$ are general divisors, satisfies the
inequality
$$
\mathop{\rm mult}\nolimits_BZ>4n^2.
$$
(The multiplicity is understood in the usual sense, see
\cite[Chapter 2]{Pukh13a}.) In fact, the assumptions can be
somewhat relaxed. The following claim is true.\vspace{0.1cm}

{\bf Proposition 2.4.} {\it Let $X$ be a variety with quadratic
singularities of rank at least 4, and assume that $\mathop{\rm
codim}(\mathop{\rm Sing}X\subset X)\geqslant 4$. Assume further
that a certain divisor $E$ over $X$ is a non canonical singularity
of the pair $(X,\frac{1}{n}\Sigma)$ with the centre
$B\subset\mathop{\rm Sing}X$, where $\Sigma$ is a mobile linear
system. Then the self-intersection $Z$ of the system $\Sigma$
satisfies the inequality}
$$
\mathop{\rm mult}\nolimits_BZ>4n^2.
$$

{\bf Proof.} We only point out what should be modified in the
arguments of \cite[Section 3]{EP}. It follows from Proposition 1.4
that the technique of counting multiplicities works without
changes under the relaxed assumptions about the rank of quadratic
singularities. Furthermore, in \cite[Section 3]{EP} it is claimed
erroneously that the Noether-Fano inequality has the form
\begin{equation}\label{18.11.2015.2}
\sum^K_{i=1}p_i\nu_i>n\left(\sum^K_{i=1}p_i\delta_i\right),
\end{equation}
where $p_i$ is the number of paths in the oriented graph of the
resolution of the singularity $E$ from the top vertex to the
vertex $i$ (the meaning of all notations is exactly the same as in
Subsection 2.1 of the present section). In fact, in the inequality
(\ref{18.11.2015.2}) instead of $p_1$ the coefficients $r_i$,
introduced in Subsection 2.1, must be used. After the replacement
of the coefficients $p_i$ by the coefficients $r_i$ all the
arguments in \cite[Section 3]{EP} work as they are and prove
Proposition 2.4.


\section{Hypersurfaces with non-isolated singularities}

In this section we prove Theorem 0.3. The procedure of estimating
the codimension of the set of hypersurfaces in the projective
space with a singular set of a positive dimension, depends on the
type of that singular set. In Subsection 3.1 we consider some
simple cases (for instance, when the singular set is a line),
where the codimension of the set of hypersurfaces with a singular
set of the given type can be directly estimated or explicitly
computed. In Subsection 3.2 we develop a technique that makes it
possible to estimate the codimension of the set of hypersurfaces
with at least finite, but sufficiently large set of singular
points. In Subsection 3.3 we apply this technique and complete the
proof of Theorem 0.3.\vspace{0.3cm}

{\bf 3.1. The sets of singular hypersurfaces.} Let ${\mathbb P}^N$
be the projective space with homogeneous coordinates
$(x_0:x_1:\dots:x_N),N\geqslant 3$, and ${\cal P}_{N,d}=
H^0({\mathbb P}^N,{\cal O}_{{\mathbb P}^N}(d))$ the linear space
of homogeneous polynomials of degree $d$. For $f\in{\cal P}_{N,d}$
the set of singular points of the hypersurface $\{f=0\}$ we denote
by the symbol $\mathop{\rm Sing}(f)$. Set
$$
{\cal P}^{(i)}_{N,d}=\{f\in{\cal P}_{N,d}\,|\,\mathop{\rm
dim}\mathop{\rm Sing}(f)\geqslant i\}.
$$
These are closed subsets of the space ${\cal P}_{N,d}$, and for
$i\geqslant j$ we have ${\cal P}^{(j)}_{N,d} \subset{\cal
P}^{(i)}_{N,d}$.\vspace{0.1cm}

{\bf Example 3.1.} Let ${\cal P}^{\rm line}_{N,d}$ be the closed
subset of the space ${\cal P}_{N,d}$, consisting of polynomials
$f$ such that the set $\mathop{\rm Sing}(f)$ contains a line in
${\mathbb P}^N$. Fixing a line $L\subset{\mathbb P}^N$, we may
assume that $L=\{x_2= \dots=x_N=0\}$, so that the condition
$L\subset\mathop{\rm Sing}(f)$ is equivalent to the set of
equalities
$$
\left.\frac{\partial f}{\partial
x_0}\right|_L\equiv\left.\frac{\partial f}{\partial x_1}\right|_L
\equiv\dots\equiv\left.\frac{\partial f}{\partial
x_N}\right|_L\equiv 0,
$$
whence, taking into account the dimension of the Grassmanian of
lines in ${\mathbb P}^N$, we obtain the equality
$$
\mathop{\rm codim}({\cal P}^{\rm line}_{N,d}\subset{\cal
P}_{N,d})=(d-2)N+3.
$$
The following claim is true, which immediately implies Theorem
0.3.\vspace{0.1cm}

{\bf Theorem 3.1.} {\it The following inequality holds:}
$$
\mathop{\rm codim}({\cal P}^{(1)}_{N,d}\subset{\cal
P}_{N,d})\geqslant (d-2)N.
$$

{\bf Remark 3.1.} It seems that the inequality of Theorem 3.1 can
be improved, replacing its right hand side by $(d-2)N+3$, after
which it would become precise. However, the proof below is
insufficient for that purpose. In any case the claim of Theorem
3.1 is much stronger than what we need in this
paper.\vspace{0.1cm}

{\bf Proof of Theorem 3.1.} Let ${\cal
P}^{(i,k)}_{N,d}\subset{\cal P}_{N,d}$ be the closure of the set,
consisting of polynomials $f$ such that $\mathop{\rm Sing}(f)$
contains an irreducible component $C$ of dimension $i\geqslant 1$,
the linear span of which $\langle C\rangle$ is a $k$-plane in
${\mathbb P}^N$, $k\geqslant i$. For instance, ${\cal
P}^{(1,1)}_{N,d}={\cal P}^{\rm line}_{N,d}$. Obviously,
$$
{\cal P}^{(i)}_{N,d}=\bigcup^N_{k=i}{\cal P}^{(i,k)}_{N,d},
$$
so that in order to estimate the codimension of the set ${\cal
P}^{(i)}_{N,d}$, it is sufficient to estimate the codimension of
each set ${\cal P}^{(i,k)}_{N,d}$, $k=i,\dots,N$. Furthermore, for
a $k$-plane $P\subset{\mathbb P}^N$ consider the set ${\cal
P}^{(i,k)}_{N,d}(P)$, which is the closure of the subset,
consisting of polynomials $f$ such that $\mathop{\rm Sing}{f}$
contains an irreducible component $C$ of dimension $i$ and
$\langle C\rangle=P$. The following fact is obvious.\vspace{0.1cm}

{\bf Proposition 3.1.} {\it The following inequality holds:}
$$
\mathop{\rm codim}\left({\cal P}^{(i,k)}_{N,d}\subset{\cal
P}_{N,d}\right)\geqslant \mathop{\rm codim}({\cal
P}^{(i,k)}_{N,d}(P)\subset{\cal P}_{N,d})-(k+1)(N-k).
$$

Finally, let ${\cal P}^{(i,k;l)}_{N,d}(P)\subset{\cal
P}^{(i,k)}_{N,d}(P)$ be the closure of the subset, consisting of
such $f$ that (in terms of the definition of the set ${\cal
P}^{(i,k)}_{N,d}(P)$) the set of singular points $\mathop{\rm
Sing}(f|_P)$ contains an irreducible component $B$ of dimension
$l$ such that $C\subset B\subset P$. In particular, $l\geqslant i$
and
$$
{\cal P}^{(i,k)}_{N,d}(P)=\bigcup^k_{l=i}{\cal
P}^{(i,k;l)}_{N,d}(P).
$$
Now by Proposition 3.1 the claim of Theorem 3.1 follows from the
system of inequalities
\begin{equation}\label{21.09.2015.1}
\mathop{\rm codim}({\cal P}^{(1,k;l)}_{N,d}(P)\subset{\cal
P}_{N,d})\geqslant(d-2)N+(k+1)(N-k),
\end{equation}
which we will prove for all $1\leqslant l\leqslant k\leqslant N$
and a fixed $k$-plane $P\subset{\mathbb P}^N$, given by the
equations $\{x_{k+1}=\dots=x_N=0\}$.\vspace{0.1cm}

{\bf Example 3.2.} Consider the case $l=k=2$. In that case $P$ is
a plane, $P\subset\{f=0\}$ and the closed set $\mathop{\rm
Sing}(f)$ contains an irreducible plane curve $C\subset P$ of
degree $q\geqslant 2$. This gives $(d+1)(d+2)/2$ independent
conditions on the coefficients of the polynomial $f|_P$ (they all
vanish) and $(N-2$) polynomials
$$
\left.\frac{\partial f}{\partial x_3}\right|_P,\dots,
\left.\frac{\partial f}{\partial x_N}\right|_P
$$
vanish on the curve $C$. Note that the coefficients of the
polynomials $f|_P$, $\partial f/\partial x_i|_P$, $i=3,\dots,N$ up
to a non-zero integral factor are distinct coefficients of the
polynomial $f$. We may assume that at least one of the polynomials
$\partial f/\partial x_i|_P$ is not identical zero, say, $\partial
f/\partial x_3|_P\not\equiv 0$. Then the curve $C$ is an
irreducible component of the plane curve $\{\partial f/\partial
x_3|_P=0$\}. Fixing the polynomial $\partial f/\partial x_3|_P$,
we finally obtain
$$
\frac{(d+1)(d+2)}{2}+(N-3)\left(qd-\frac{q(q-1)}{2}\right)
$$
independent conditions on the coefficients of the polynomial $f$,
where $2\leqslant q\leqslant d-1$. It ie easy to see that this
number satisfies the inequality
(\ref{21.09.2015.1}).\vspace{0.1cm}

{\bf Example 3.3.} Consider the case $l=1,k=2$. In that case
$\mathop{\rm Sing}(f)$ contains an irreducible plane curve
$C\subset P$ of degree $q\geqslant 2$, but $f|_P\not\equiv 0$, so
that $\{f|_P=0\}$ is a reducible plane curve of degree $d$,
containing $C$ as a double component, so that $2q\leqslant d$. An
easy dimension count gives
$$
\frac12(5q^2-(4d+3)q+d^2+3d+4)
$$
independent conditions on the coefficients of the polynomial
$f|_P$. The minimum of the last expression is attained for $q=2$.
Now the fact that the polynomials $\partial f/\partial x_i|_P$,
$i=3,\dots,N$, vanish on the curve $C$, gives in addition at least
$(N-2)(2d+1)$ independent conditions on the coefficients of $f$.
As a result, we get
$$
\mathop{\rm codim}\left({\cal
P}^{(1,2;1)}_{d,N}(P)\right)\geqslant\frac{d(d-5)}{2}+(N-2)(2d+1)+9,
$$
and it is easy to check that the inequality (\ref{21.09.2015.1})
is satisfied.\vspace{0.1cm}

Starting from this moment, we assume that $k=3$. Recall the
following\vspace{0.1cm}

{\bf Definition 3.1.} (See \cite[Section 3]{Pukh01} or
\cite[Chapter 3]{Pukh13a}). A sequence of homogeneous polynomials
$g_1,\dots,g_m$ of arbitrary degrees on the projective space
${\mathbb P}^e$, $e\geqslant m+1$, is called a {\it good
sequence}, and an irreducible subvariety $W\subset{\mathbb P}^e$
of codimension $m$ is its {\it associated subvariety}, if there
exists a sequence of irreducible subvarieties $W_j\subset{\mathbb
P}^e$, $\mathop{\rm codim}W_j=j$ (in particular, $W_0={\mathbb
P}^e$), such that:
\begin{itemize}

\item $g_{j+1}|_{W_j}\not\equiv 0$ for $j=0,\dots,m+1$,

\item $W_{j+1}$ is an irreducible component of the closed
algebraic set $g_{j+1}|_{W_j}=0$,

\item $W_m=W$.

\end{itemize}

\noindent A good sequence can have more than one associated
subvarieties, but their number is bounded from above by a constant
depending on the degrees of the polynomials $q_j$ only (see
\cite[Section 3]{Pukh01}).\vspace{0.1cm}

Let us consider two more examples, similar to Examples 3.2 and
3.3.\vspace{0.1cm}

{\bf Example 3.4.} Let us consider the case $l=k$. This case
generalizes Example 3.2. We have $f|_P\equiv 0$, which gives
${k+d\choose d}$ independent conditions on the coefficients of
$f$. Since the polynomials
$$
\left.\frac{\partial f}{\partial x_{k+1}}\right|_P,\dots,
\left.\frac{\partial f}{\partial x_N}\right|_P
$$
vanish identically on $C$ and the curve $C$ is an irreducible
component of the set $\mathop{\rm Sing}(f)$, from those
polynomials we can choose $(k-1)$ ones that form a good sequence
with the curve $C$ as an associated subvariety (in particular,
$N-k\geqslant k-1$). Fixing these polynomials, for each of the
remaining $(N+1-2k)$ polynomials we get the condition
\begin{equation}\label{22.09.2015.1}
\left.\left(\left.\frac{\partial f}{\partial
x_i}\right|_P\right)\right|_C\equiv 0,
\end{equation}
where the curve $C$, as one of the associted subvarieties of the
fixed good sequence, can be assumed to be fixed. In \cite[Section
3]{Pukh01} it was shown that the condition (\ref{22.09.2015.1})
defines a closed subset of codimension at least $(d-1)k+1$.
Therefore,
$$
\mathop{\rm codim}({\cal P}^{(i,k;k)}_{N,d}\subset{\cal P}_{N,d})
\geqslant {k+d\choose d}+(N+1-2k)((d-1)k+1),
$$
and elementary computations show that the inequality
(\ref{21.09.2015.1}) holds.\vspace{0.1cm}

{\bf Example 3.5.} Let us consider the case $l=k-1$. This case
generalizes Example 3.3. Here the hypersurface $\{f|_P=0\}$ has a
multiple irreducible non-degenerate component of degree $q$, where
$2q\leqslant d$, so that the coefficients of the polynomial $f|_P$
belong to a closed subset of codimension
$$
{k+d\choose k}-{k+d-2q\choose k}-{k+q\choose k}
$$
in the space ${\cal P}_{d,k}$. Furthermore, since the curve $C$ is
an irreducible component of the set $\mathop{\rm Sing}(f)$, from
the set of polynomials
$$
\left.f|_P, \frac{\partial f}{\partial x_{k+1}}\right|_P,\dots,
\left.\frac{\partial f}{\partial x_N}\right|_P
$$
we can choose a good sequence, starting with $f|_P$, for which the
curve $C$ will be an associated subvariety. In particular, the
estimate $N+2\geqslant 2k$ holds. Fixing the polynomials of that
sequence, we may assume the curve $C$ to be fixed. Now we argue as
in Example 3.4 and obtain, in addition to the conditions on the
coefficients of the polynomial $f|_P$, also $(N+2-2k)((d-1)k+1)$
more independent conditions on the coefficients of the polynomial
$f$. An elementary, although tedious, check shows that the
inequality (\ref{21.09.2015.1}) is satisfied.\vspace{0.1cm}

In order to prove the inequality (\ref{21.09.2015.1}) in the case
$l\leqslant k-2$, we need a new technique, which is developed
below.\vspace{0.3cm}


{\bf 3.2. Linearly independent points.} The following claim is
true.\vspace{0.1cm}

{\bf Lemma 3.1.} {\it Assume that $d\geqslant 3$. For any set of
$m$ linearly independent points $p_1,\dots,p_m\in{\mathbb P}^N$,
$m\leqslant N+1$, the condition
$$
\{p_1,\dots,p_m\}\subset\mathop{\rm Sing}(g),
$$
$g\in{\cal P}_{N,d}$, defines a linear subspace of codimension
$m(N+1$) in ${\cal P}_{N,d}$.}\vspace{0.1cm}

{\bf Proof.} We may assume that
$$
p_1=(1:0:0\dots:0),\quad p_2=(0:1:0:\dots:0)
$$
etc. correspond to the first $m$ vectors of the standard basis of
the linear space ${\mathbb C}^{N+1}$. The condition
$p_i\in\mathop{\rm Sing}(g)$ means vanishing the coefficients at
the monomials $x^d_{i-1},x^{d-1}_{i-1}x_j$, for all $j\neq i-1$.
For $d\geqslant 3$ all those $m(N+1)$ monomials are distinct.
Q.E.D. for the lemma.\vspace{0.1cm}

Now let us consider an arbitrary linear subspace
$\Pi\subset{\mathbb P}^N$ of codimension $r+1$, where $r\geqslant
1$, given by a system of $r+1$ equaltions
$$
l_0(x)=0,\,\,l_1(x)=0,\dots,l_r(x)=0,
$$
where $l_0,\dots,l_r$ are linearly independent forms. For every
$i=1,\dots,r$ let us fix an arbitrary tuple of distinct constants
$\lambda_{i0},\dots,\lambda_{i,d-1}\in{\mathbb C}$; we assume that
$\lambda_{i0}=0$ for all  $i=1,\dots,r$. Now for any integer point
$$
\underline{e}=(e_1,\dots,e_r)\in{\mathbb Z}^r_+,\quad e_i\leqslant
d-1,
$$
by the symbol $\Theta(\underline{e})$ we denote the linear
subspace
$$
\{l_i(x)-\lambda_{i,e_i}l_0(x)=0\,|\,i=1,\dots,r\}\subset{\mathbb
P}^N
$$
of codimension $r$. Obviously, $\Theta(\underline{e})\supset\Pi$.
Set
$$
|\underline{e}|=e_1+\dots+e_r\in{\mathbb Z}_+.
$$
For each tuple $\underline{e}\in{\mathbb Z}^r_+$ with
$|\underline{e}|\leqslant d-3$ consider an arbitrary set
$$
S(\underline{e})=\{p_1(\underline{e}),\dots,p_m(\underline{e})\}
\subset\Theta(\underline{e})\backslash\Pi
$$
of $m$ linearly independent points (so that $m\leqslant
N-r+1$).\vspace{0.1cm}

{\bf Proposition 3.2.} {\it The set of conditions
$$
S(\underline{e})\subset\mathop{\rm
Sing(g|_{\Theta(\underline{e})})},
$$
$\underline{e}\in{\mathbb Z}^r_+$, $|\underline{e}|\leqslant d-3$,
defines a linear subspace of codimension
$$
m(N-r+1)|\Delta|
$$
in ${\cal P}_{N,d}$, where
$$
\Delta=\{e_1\geqslant 0,\dots,e_r\geqslant
0,e_1+\dots+e_r\leqslant d-3\}\subset{\mathbb R}^r
$$
is an integral simplex and $|\Delta|$ means the number of integral
points in that simplex, $|\Delta|=\sharp(\Delta\cap{\mathbb
Z}^r)$.}\vspace{0.1cm}

{\bf Proof.} We may assume that $l_0=x_0$,
$l_1=x_1,\dots,l_r=x_r$. In order to simplify the formulas, we
will prove the affine version of the proposition: set
$v_1=x_1/x_0,\dots, v_r=x_r/x_0$ and $u_i=x_{r+i}/x_0$,
$i=1,\dots,N-r$. In the affine space ${\mathbb A}^N\subset{\mathbb
P}^N$, ${\mathbb A}^N={\mathbb P}^N\backslash\{x_0=0\}$ with
coordinates $(u,v)=(u_1,\dots,u_{N-r}$, $v_1,\dots,v_r$) the
affine spaces
$A(\underline{e})=\Theta(\underline{e})\backslash\Pi$ are
contained entirely:
$$
A(\underline{e})=\Theta(\underline{e})\cap{\mathbb A}^N,
$$
so that $S(\underline{e})\subset A(\underline{e})$ for all
$\underline{e}$. Obviously,
$$
A(\underline{e})=
\{v_1=\lambda_{1,e_1},\dots,v_r=\lambda_{r,e_r}\}\subset{\mathbb
A}^N
$$
is a $(N-r)$-plane, which is parallel to the coordinate
$(N-r)$-plane $(u_1,\dots,u_{N-r},0,\dots,0)$. Let us write the
polynomial $g$ in terms of the affine coordinates $(u,v)$ in the
following way:
$$
g(u,v)=\sum_{\underline{e}\in{\mathbb
Z}^r_+,|\underline{e}|\leqslant d}
g_{e_1,\dots,e_r}(u)\prod^r_{i=1}\prod^{e_i-1}_{j=0}(v_i-\lambda_{ij})
$$
(if $e_i=0$, then the corresponding product is assumed to be equal
to  1). Here $g_{\underline{e}}(u)=g_{e_1,\dots,e_r}(u)$ is an
affine polynomial in  $u_1,\dots,u_{N-r}$ of degree $\mathop{\rm
deg}g_{\underline{e}}\leqslant d-|e|$. For the fixed
$\lambda_{ij}$ this presentation is unique. By Lemma 3.1, the
condition
$$
S(\underline{0})=S(0,\dots,0)\subset\mathop{\rm
Sing}(g|_{A(\underline{0})})
$$
defines a linear subspace of codimension $m(N-r+1)$ in the space
of polynomials ${\cal P}_{N-r,d}$. However, it is easy to see that
$$
g|_{A(\underline{0})}=g_{0,\dots,0}(u),
$$
since for $\underline{e}\neq\underline{0}$ in the product
$$
\prod^r_{i=1}\prod^{e_i-1}_{j=0}(v_i-\lambda_{ij})
$$
there is at least one factor $(v_i-\lambda_{i0})=v_i$, which
vanishes when restricted onto the $(N-r)$-plane
$A(\underline{0})$. Therefore, the condition
$S(\underline{0})\mathop{\rm Sing}(g|_{A(\underline{0})})$ imposes
on the coefficients of the polynomial $g_{0,\dots,0}(u)$ precisely
$m(N-r+1)$ independent linear conditions, whereas the polynomials
$g_{\underline{e}}(u)$ for $\underline{e}\neq\underline{0}$ can be
arbitrary.\vspace{0.1cm}

Now let us complete the proof of Proposition 3.2 by induction on
$|\underline{e}|$. More precisely, for any $a\in{\mathbb Z}_+$ set
$$
\Delta_a=\{e_1\geqslant 0,\dots,e_r\geqslant
0,\,\,e_1+\dots+e_r\leqslant a\}\subset{\mathbb R}^r,
$$
so that $\Delta=\Delta_{d-3}$, and let us prove the claim of
Proposition 3.2 in the following form: for every
$a=0,\dots,d-3$\vspace{0.1cm}

\parshape=1
2cm 11cm \noindent $(*)_a$ the set of conditions
$$
S(\underline{e})\subset\mathop{\rm
Sing}(g|_{\Theta(\underline{e})}),
$$
$\underline{e}\in{\mathbb Z}^r_+$, $|\underline{e}|\leqslant a$,
defines a linear subspace of codimension $m(N-r+1)|\Delta_a|$ in
${\cal P}_{N,d}$, where the restrictions are imposed on the
coefficients of the polynomials $g_{\underline{e}}(u)$ for
$\underline{e}\in\Delta_a$, whereas for
$\underline{e}\not\in\Delta_a$ the polynomials
$g_{\underline{e}}(u)$ can be arbitrary.\vspace{0.1cm}

The case $a=0$ has already been considered, so we assume that
$a\leqslant d-4$ and the claims $(*)_j$ have been shown for
$j=0,\dots,a$. Let us show the claim $({*})_{a+1}$. Let
$\underline{e}\in{\mathbb Z}^r_+$ be an arbitrary multi-index,
$|\underline{e}|=a+1$. The restriction onto the affine subspace
$A(\underline{e})$ means the substitution
$v_1=\lambda_{1,e_1},\dots$, $v_r=\lambda_{r,e_r}$. For that
reason the polynomial $g_{\underline{e}(u)}$ comes into the
restriction $g|_{A(\underline{e})}$ with a non-zero coefficient
$$
\alpha_e=\prod^r_{i=1}\prod^{e_i-1}_{j=0}(\lambda_{i,e_i}-\lambda_{ij}).
$$
On the other hand, for $\underline{e}'\neq\underline{e}$,
$|\underline{e}'|\geqslant a+1$ the product
$$
\prod^r_{i=1}\prod^{e'_i-1}_{j=0}(\lambda_{i,e_i}-\lambda_{ij}).
$$
is equal to zero, as for at least one index $i\in\{1,\dots,r\}$ we
have $e'_i> e_i$ and therefore that product contains a zero
factor. So $g|_{A(\underline{e})}$ is the sum of the polynomial
$\alpha_eg_{\underline{e}}$ and a linear combination of the
polynomials $g_{\underline{e}'}$ with $|\underline{e}'|\leqslant
a$ with constant coefficients. Now, fixing the polynomials
$g_{\underline{e}'}$ с $|\underline{e}'|\leqslant a$, we see that
the condition
$$
S(\underline{e})\subset\mathop{\rm Sing}(g|_{A(\underline{e})})
$$
defines an {\it affine} (generally speaking, not a linear)
subspace of codimension $m(N-r+1)$ of the space of polynomials
$g_{\underline{e}}(u_1,\dots,u_{N-r})$ of degree at most $d-|e|$,
the corresponding linear space of which is given by the condition
$$
S(\underline{e})\subset\mathop{\rm Sing}g_{\underline{e}}(u).
$$
Note that on the coefficients of other polynomials
$g_{\underline{e'}}$ with $|\underline{e}'|=a+1$ no restrictions
are imposed.\vspace{0.1cm}

This completes the proof of the claim $(*)_a$ for all
$a=0,\dots,d-3$. Q.E.D. for Proposition 3.2.\vspace{0.3cm}


{\bf 3.3. End of the proof of Theorem 3.1.} Let
$$
\Theta=\Theta[l_0,\dots,l_r;\lambda_{i,j},i=1,\dots,r,j=0,\dots,d-1]=
\{\Theta(\underline{e})\,|\,\underline{e}\in\Delta\}
$$
be a set of linear subspaces of codimension $r$ in ${\mathbb
P}^N$, considered in Proposition 3.2. We define the subset
$$
{\cal P}_{N,d}(\Theta)\subset{\cal P}_{N,d}
$$
by the following condition: for every subspace
$\Theta(\underline{e})$ with $|\underline{e}|\leqslant d-3$ there
is a set
$S(\underline{e})\subset\Theta(\underline{e})\backslash\Pi$,
consisting of $m$ linearly independent points, such that
$S(\underline{e})\subset\mathop{\rm
Sing}(g|_{\Theta(\underline{e})})$.\vspace{0.1cm}

{\bf Proposition 3.3.} {\it The following inequality is true:}
$$
\mathop{\rm codim}({\cal P}_{N,d}(\Theta)\subset{\cal
P}_{N,d})\geqslant m|\Delta|.
$$

{\bf Proof} is obtained by the obvious dimension count: the
subspaces $\Theta(\underline{e})$ are fixed, so that every point
$p_i(\underline{e})$ varies in a $(N-r)$-dimensional family.
Q.E.D. for the proposition.\vspace{0.1cm}

Let us complete, finally, the proof of Theorem 3.1. Consider the
set ${\cal P}^{(1,k;l)}_{N,d}(P)$, where $P\subset{\mathbb P}^N$
is a fixed $k$-plane $\{x_{k+1}=\dots=x_N=0\}$, and $l\leqslant
k-2$. We apply proposition 3.3 to the space $P$ instead of
${\mathbb P}^N$ and to the space of polynomials ${\cal P}_{k,d}$
instead of ${\cal P}_{N,d}$. For an arbitrary set
$\Theta=\{\Theta(\underline{e})\,|\, \underline{e}\in\Delta\}$ of
linear subspaces of codimension $l$ in $P={\mathbb P}^k$ let
$$
{\cal P}^{(1,k;l)}_{N,d}(P,\Theta)\subset{\cal P}^{(1,k;l)}_{N,d}
$$
be the set of polynomials $f\in{\cal P}^{(1,k;l)}_{N,d}$ such that
the set $\mathop{\rm Sing}(f|_P)$ has an irreducible component $Q$
of dimension $l$, containing a curve $C\subset\mathop{\rm
Sing}(f)$, and such that it is in general position with the
subspaces from the set $\Theta$: for all $\underline{e}\in\Delta$
the set $\Theta(\underline{e})\cap Q$ contains $(k-l+1)$ linearly
independent points. Since $\langle Q\rangle=\langle C\rangle=P$,
the subset ${\cal P}^{(1,k;l)}_{N,d}(P,\Theta)$ is a Zariski open
subset of the set ${\cal P}^{(1,k;l)}_{N,d}$, so that the
inequality (\ref{21.09.2015.1}) will be shown is we show it for
${\cal P}^{(1,k;l)}_{N,d}(P,\Theta)$ instead of ${\cal
P}^{(1,k;l)}_{N,d}$. By Proposition 3.3, applied to the space $P$,
the condition $f\in{\cal P}^{(1,k;l)}_{N,d}(P,\Theta)$ imposes on
the coefficients of the polynomial $f|_P$ at least
$(k-l+1)|\Delta|$ independent conditions. Furthermore, from the
set of $(N+1$) polynomials
$$
\left.\frac{\partial f}{\partial x_0}\right|_P,\dots,
\left.\frac{\partial f}{\partial x_N}\right|_P
$$
we may select a good sequence of $(k-1)$ polynomials, with a
certain curve $C$, $\langle C\rangle=P$, as an associated
subvariety, and moreover, this can be done in such a way that the
first $(k-l)$ polynomials in that sequence are chosen among the
polynomials
$$
\left.\frac{\partial f}{\partial x_0}\right|_P,\dots,
\left.\frac{\partial f}{\partial x_k}\right|_P
$$
(and some subvariety $Q\supset C$, $Q\subset P$ of dimension $l$
is an associated subvariety of that subsequence), whereas the
following $(l-1)$ polynomials are chosen among the polynomials
$$
\left.\frac{\partial f}{\partial x_{k+1}}\right|_P,\dots,
\left.\frac{\partial f}{\partial x_N}\right|_P.
$$
Fixing the polynomial $f|_P$ and the other polynomials of the good
sequence, we may assume the curve $C\subset\mathop{\rm Sing}(f)$
of singular points to be fixed. Now the condition $\partial
f/\partial x_i|_C\equiv 0$ for every $i\in\{k+1,\dots,N\}$, which
did not get into the good sequence, give in addition
$(N+1-k-l)((d-1)k+1)$ independent conditions on the coefficients
of the polynomial $f$. An elementary, although tedious, check
shows that the inequality
$$
(k-l+1)|\Delta|+(N+1-k-l)((d-1)k+1)\geqslant(d-2)N+(k+1)(N-k)
$$
holds for all the values $k,l$ under consideration, which
completes the proof of the inequality (\ref{21.09.2015.1}) and of
Theorem 3.1, and therefore, of Theorem 0.3.\vspace{0.1cm}

{\bf Remark 3.2.} It is easy to see that the worst estimate for
the codimension of ${\cal P}^{(1,k;l)}_{N,d}(P)$ corresponds to
the case $k=N$ and $l=1$, that is, the hypersurface $\{f=0\}$ has
a non-degenerate curve of singular points. In that case
Proposition 3.3 yields the inequality
$$
\mathop{\rm codim}({\cal P}^{(1,k;1)}_{N,d}\subset{\cal
P}_{N,d})\geqslant(d-2)N.
$$
It seems hardly probable that the presence of a non-degenerate
curve of singular points imposes on the coefficients of the
polynomial $f$ less (although slightly less) independent
conditions than the presence of a line consisting of singular
points (when the estimate for the codimension is precise). And
indeed, when we apply Proposition 3.3, we essentially replace a
curve, consisting of singular points, by a finite set of singular
points (although it is quite a large set). Probably, the technique
used in the proof of Theorem 3.1 can be improved and for the case
of a non-degenerate curve of singular points a more precise
estimate could be obtained. This is what was meant in Remark 3.1.

\begin{flushleft}
Department of Mathematical Sciences,\\
The University of Liverpool
\end{flushleft}

\noindent{\it pukh@liverpool.ac.uk}

\end{document}